\documentclass[a4paper,12pt]{article}
\usepackage{amssymb}
\usepackage[french,english,german]{babel}
\def \d {{\mathrm {d}}}

\def \b {{\bigskip}}
\def \m {{\medskip}}

\title{From Schouten to Mackenzie: notes on brackets\footnote{To appear in the {\it Journal of Geometric Mechanics}, 2021.}}
             
\author{Yvette Kosmann-Schwarzbach}

\date{}

\begin{document}

\maketitle
{\hfill {{\it In Memory of Kirill Mackenzie (1951--2020)}}}

\b

\b
\begin{small}
\noindent {\it Abstract}  
In this paper, dedicated to the memory of Kirill Mackenzie, I relate the origins and early development of the theory of graded Lie brackets, first in the publications on differential geometry of Schouten, Nijenhuis, and Fr\"olicher--Nijenhuis, then in the work of Gerstenhaber and Nijenhuis--Richardson in cohomology theory. 

\b

\b

\end{small}

\b 
\section*{Introduction}
Many aspects of the contributions of Kirill Mackenzie to the theory of Lie groupoids and Lie algebroids, in particular to that of the double and multiple structures, are already classical and often cited.   
As a tribute to his memory, I shall survey the history of the emergence of some theories that were tools -- among many others -- in the development and description of his mathematics.

In his {\it General Theory of Lie Groupoids and Lie Algebroids}, in 2005, Mackenzie inserted a section of ``Notes'' at the end of each chapter.
Yet he advised his readers that his notes were ``not comprehensive historical surveys''. After the chapter ``Cohomology and Schouten calculus for Lie algebroids'', he wrote: 
\begin{quotation}
\noindent ``The Schouten bracket of multivector fields on a manifold goes back  to
 Schouten [1954] and Nijenhuis [1955].
The same algebraic structures on the cohomology of associative algebras were introduced by Gerstenhaber [1963].''
\end{quotation} 
I shall try to expand on these elements of a story, and survey some of the early developments in the theory of the various graded brackets that Mackenzie cited and used\footnote{The further developments of these topics are, of course, of great interest, but I have refrained from presenting them here, as I also stopped short of dealing with the history of the Whitehead bracket (for which see J. D. Stasheff, ``Brackets by any other name'', in this volume), since it does not appear in Mackenzie's work.}.


The role of the Schouten bracket is now essential in the geometry of manifolds as well as in quantization theory\footnote{When Maxim Kontsevitch was named professor at the Institut des Hautes \'Etudes Scientifiques in 1995, he deemed it necessary to spend the first five minutes of his inaugural lecture delivering an overview of the Schouten bracket before he went on to report on his latest work.}. I shall outline its introduction in the work of the Dutch mathematician Jan Schouten  and analyse the subsequent contribution of his student Albert Nijenhuis, so that this bracket can rightfully be called ``the Schouten--Nijenhuis bracket'' while it can also be 
called ``the Gerstenhaber bracket''. For good measure, I shall sketch the origins of the Fr\"olicher--Nijenhuis and Nijenhuis--Richardson brackets.

\m
Below, I often write ``tensor'' for ``tensor field''. The term ``concomitant'' or ``comitant'' was used in the early publications for what is now called a bracket. It was still in use in some papers on geometry in the 1970's, while the use of ``concomitant'' became frequent in statistics, though with a different meaning.
The use of brackets $[~,~]$ or parentheses $(\,,\,)$ to denote the Poisson brackets in mechanics goes back to their inventors Joseph Louis Lagrange and Sim\'eon Denis Poisson, but the notations are frequently exchanged in their early writings of 1810 and 1811. These notations and the braces $\{\,,\, \}$ have all been used at various times and in various contexts for the other pairings that were later introduced in geometry and algebra, such as those whose invention I shall trace here.

\subsection*{Dramatis personae}
\begin{quotation}
\noindent ``The first basic example of graded Lie algebras was provided by Nijenhuis (1955) and then by Fr\"olicher and Nijenhuis (1956) [...]. The basic role that this object plays in the theory of deformation of algebraic structures was discovered by Gerstenhaber (1963, 1964) in a fundamental series of papers, while Spencer and collaborators were developing applications to pseudogroup structures on manifolds [...].
The general role of graded Lie algebras in  deformation theory was presented by Nijenhuis and Richardson (1964).''
\end{quotation}
%
So wrote Lawrence Corwin, Yuval Ne'eman and Shlomo Sternberg in 1975 as they introduced their up-to-date survey on ``Graded Lie algebras in mathematics and physics'' \cite{CNS1975}, not only reviewing ``the recent applications of supersymmetry in the physics of particles and fields'' but also presenting ``a short survey of the role played by graded Lie algebras in mathematics''. I intend to follow the development of this subject in the contributions of Schouten and of the other mathematicians whom they identified.

\subsection*{``A master at tensors''}
Jan Arnoldus Schouten  (1883--1971) began his studies in electrotechnical engineering in Delft in 1901, and only turned to mathematics in 1912. In his dissertation, ``Grundlagen der Vektor- und Affinoranalysis'' (Foundations of vector and affinor analysis) (1914) \cite{Schouten1914}, he applied the classification methods, inspired by Klein's Erlangen program, to the study of  ``direct quantities'', the precursors of tensors, elucidating the confusing terminology and notation then in use. At the time, tensors were called ``affinors'',  while the term ``tensors'' was reserved for those affinors that are symmetric. It was Einstein who adopted the term tensor in its present sense and the usage changed slowly thereafter. 

Schouten served as professor in Delft University from 1914 until 1943\footnote{Schouten's first assistant in Delft was a woman, Johanna Manders. Her doctorate in 1919 may have been the first doctor's degree in mathematics conferred on a woman in Delft or even in the Netherlands.}. In 1948 he accepted a professorship in Amsterdam, where he co-founded the Mathematische Centrum, became its director and remained a member until 1968. 
The website zbMATH lists 204 publications, including 23 books, bearing his name as author that were  
published between 1914 and 1978. Part of Schouten's 
{\oe}uvre was surveyed in the  article ``Schouten and the tensor calculus'' \cite{S1978} by Dirk Struik (1894-2000) who had been his assistant in Delft from 1917 to 1923. ``A master at tensors'' is the title of the obituary for Schouten \cite{N1972} that Nijenhuis wrote in 1972.

\subsection*{A gifted doctoral student}
Albert Nijenhuis  (1926-2015) was  one of Schouten's doctoral students in Amsterdam where he completed his thesis in 1952, ``Theory of the geometric object'' \cite{N1952}.
He left for the United States where he held positions at the Institute for Advanced Studies in Princeton, then at the University of Chicago, and in 1957 he was assistant professor of mathematics at the University of Washington in Seattle. He returned to Europe  to take up a position at the University of Amsterdam for the year 1963-1964 and was then appointed to the University of Pennsylvania in Philadelphia. There, his research interest strayed from differential geometry and he is well known as the co-author of the book {\it Combinatorial Algorithms} which appeared in 1975. After he retired in 1987, he returned to Seattle where he remained active as affiliate professor in the mathematics department.  

Below I shall discuss in some detail the insightful paper which he wrote in 1951 before defending his thesis, as well as his 1955 paper in which he formulated in modern terms a result that figured in Schouten's 1940 paper.


\section{The Schouten--Nijenhuis bracket}
\subsection*{Schouten's concomitant, 1940}

In 1940, Schouten, in a 4-page article, \selectlanguage{german}``Ueber Differentialkomitanten zweier kontravarianter 
Gr\"ossen''\selectlanguage{english} (On differential comitants of two contravariant quantities) \cite{Schouten1940}, defined a generalization of the Lie bracket, a differential concomitant on the space of all contravariant affinors, i.e., tensors, on a manifold. 

Schouten considered the Lie derivative operator, as it was then known, acting on contravariant tensors of arbitrary valence, for which he cited  W{\l}adys{\l}aw {\'S}lebodzi{\'n}ski \cite{WS1931}, as well as the book that he had written with Struik \cite{SchoutenStruik}.
His seminal idea -- which he 
explained later \cite{Schouten2} -- was that, in the expression of the Lie derivative of a contravariant tensor field $t$ by a vector field $X$, one could look at the result as either an operator  on $t$ depending on $X$, or as an operator on $X$ depending on $t$. This idea led him to the search for a differential concomitant of a pair of geometric objects that were more general than ``vector-tensor'' pairs, namely ``tensor-tensor'' pairs. 
He found an expression in local coordinates for a differential concomitant or, in modern terms, a bracket defined for a pair of contravariant tensors, both in terms of partial derivatives and in terms of covariant derivatives with respect to a connection. He thus obtained, for two contravariant tensors of valence $a+1$ and $b+1$, respectively, the components of a contravariant tensorial quantity of valence $a+b+1$, and he showed that it indeed reduces to the Lie derivative when the first tensor is a vector. 

Then Schouten considered the case where both tensors are symmetric (respectively, anti-symmetric), and gave simplified formulas for the results, observing that the resulting tensor is symmetric (respectively, anti-symmetric). Thus were obtained both a concomitant on symmetric contravariant tensors and on skew-symmetric contravariant tensors. The restriction of Schouten's concomitant to the space of skew-symmetric contravariant tensors is the present-day ``Schouten bracket'', also called, for reasons to be explained below, the ``Schouten--Nijenhuis bracket'', of multivectors\footnote{Multivectors are called polyvectors by our Russian colleagues and others, despite the mixture of Greek and Latin roots.}.

\subsection*{An Italian mathematical extravaganza, 1953}
Differential concomitants were the main subject of Schouten's lecture, thirteen years after his short paper mentioned above. It was presented at the first International Conference on differential geometry that was organized in Italy in September 1953 by the Unione Matematica Italiana with the participation of the Consiglio Nazionale delle Ricerche, once the trauma of
the Second World War had passed. This was the first congress to be sponsored by the International Mathematical Union, and it was  a truly international gathering, with a stunning roster of speakers. Besides Schouten, almost all the well-known differential geometers from Italy, France, England, Belgium, Germany, Russia and Denmark attended.
Among them were {\v C}ech, Finsler, Haantjes\footnote{Not so well known is Johannes Haantjes (1909--1956) who had been Schouten's assistant in Delft. At the time of this conference, he was teaching at the Vrije Universiteit in Amsterdam and was a member of the Royal Dutch Academy of Sciences. His work on brackets was nearly forgotten after his death but has now become the subject of renewed interest.},
Hodge, Hopf, Kuiper, P. Libermann, Lichnerowicz, Reeb, Segre, Severi, Willmore, and Yano. One can only imagine their exchanges in the setting of this memorable gathering!
The conference was preceded by a reception in Venice, and the opening session took place the next day in Padua. Scientific sessions were held on three successive days in Venice, Bologna and Pisa, and on each day a tribute was paid, on the first day to Gregorio Ricci-Curbastro (1853--1925), honoring the memory of Tullio Levi-Civita (1873--1941) as well, on the second day to Luigi Cremona (1830--1903), and on the third day to Luigi Bianchi (1856--1928), while in the foreword to the volume of proceedings, Enrico Bompiani (1889--1975), president of the Union and honorary delegate of the IMU, stressed that ``the same tribute [should be paid] to the memory of \'Elie Cartan'' who had died two years earlier, and he also regretted the absence ``for personal reasons'' of Albert Einstein (1879--1955) and Hermann Weyl (1885--1955). The conference ended with a ceremony in Pisa. 

The proceedings appeared in 1954, with Schouten's article on concomitants, ``On the differential operators of first order in tensor calculus'' \cite{Schouten1954}, printed as a preamble to all the other communications. He recalled that ``in 1940 the present author succeeded in generalizing Lie's operator by forming a differential concomitant of two arbitrary {\it contravariant} quantities'', and he disclosed the method he had used to derive his concomitant: he discovered ``the Schouten bracket'' of contravariant tensors by requiring that it be a derivation in each argument, which is in fact the current definition of this 
concept. In this short paper, he went on to cite a paper of Emmy Noether on differential concomitants \cite{N1918} contrasting her results with the new concomitant that Nijenhuis had found in 1951, of which more below. In the last lines of his paper, 
Schouten gave the latest news about ``another generalization,'' 
conveying his excitement at the news of his former student's discovery:
\begin{quotation}
Mr. Nijenhuis informed me in his letter of 30-9-53 that he was already in the possession of another generalization. From two quantities with one upper index and alternating in the lower indices another quantity of the same kind can be formed and for these concomitants a kind of Jacobi identity exists. 
At the moment I cant [{\it sic}] give any more details.
\end{quotation}
This tantalizing last sentence of his article announced Nijenhuis's results concerning the brackets of vector-valued differential forms of arbitrary degree, results that he would publish in the second part of his article which would appear in 1955 \cite{N1955} and which he went on to develop in his joint article with Fr\"olicher in 1956 \cite{FN1956}.

\subsection*{The Schouten--Nijenhuis bracket of multivectors, 1955}
The first part of Nijenhuis's article of 1955, entitled ``Jacobi-type identities for bilinear differential concomitants of certain tensor fields'' \cite{N1955} was  devoted to the properties of the bracket of contravariant tensors which Schouten had defined in 1940. Nijenhuis showed that this bracket decomposed as the sum of two brackets, one depending only on the symmetric part of each entry and the other depending only upon the antisymmetric part. Then he
observed that, when writing ${\mathcal L}_T u$ as $[u,T]$ for the Lie derivative of a tensor field $T$ by a vector field $u$, the properties of the Lie derivation yield what he called ``a Jacobi-type identity'',
$$[u,[v,T]]- [v,[u,T]] - [[u,v],T] =0,$$
when $v$ is a vector field. 
This observation was preliminary to his main result (Theorem, p. 393) in which he proved both the {\it Jacobi identity} for the restriction of this bracket to fields of {\it symmetric} contravariant tensors,
$$[P,[Q,R]] +  [R,[P,Q]] + [Q,[R,P]] =0,$$
and a {\it graded Jacobi identity} for the restriction of this bracket to fields of {\it antisymmetric} contravariant tensors, 
$$(-1)^{r(p+1)}[P,[Q,R]] + (-1)^{q(r+1)} [R,[P,Q]] + (-1)^{p(q+1)}[Q,[R,P]] =0,$$
for $P$ of degree $p$, $Q$ of degree $q$, and $R$ of degree $r$.

In a note added in proof, Nijenhuis showed that the Schouten bracket of {\it symmetric} contravariant tensors on an analytic manifold is nothing other than the Poisson bracket of functions on its cotangent bundle.

\subsection*{Reviews and reviewers}
In {\it Mathematical Reviews} in 1951, Nijenhuis, who had not yet defended his doctoral thesis,  reviewed Schouten's ``Sur les tenseurs...'' \cite{Schouten1951} while Schouten reviewed Nijenhuis's paper \cite{N1951} that appeared shortly afterwards.

In his book, {\it Ricci-Calculus, An Introduction  to  Tensor  Analysis  and its  Geometrical  Applications} \cite{Schouten2}, published in 1954 -- a completely new version of his {\it Der Ricci-Kalk\"ul} of 1924 -- Schouten spelled out Nijenhuis's formula for the torsion of a mixed 2-tensor from \cite{N1951} 
and explained the role of the new concomitant 
$H^{\cdot \cdot \kappa}_{\lambda \mu}$
in the problem of the integrability of the field of eigenspaces
of mixed tensors, stressing that  
``Nijenhuis gives a necessary and sufficient condition for $n-1$ covariant eigenvectors of a tensor with $n$ different eigenvalues on a manifold of dimension $n$ to define a submanifold of dimension $n-1$''.
 The reviewer for Schouten's  {\it Ricci-Calculus} for the {\it Bulletin of the American Mathematical Society} was Kentaro Yano (1912--1993) who had already published many often cited books and articles on differential geometry. He concluded his entirely favorable review with: ``This  excellent  book [...] will give all the information  in the small necessary  for the development  of  differential  geometry  in the large.''

As early as 1955, Nijenhuis's name appeared in the title of an article, 
 ``On the geometrical meaning of the vanishing of the Nijenhuis tensor in an $X_{2n}$ with an almost complex structure'', which  Schouten published in collaboration with Yano \cite{SchoutenYano1955b}. It was reviewed by A. G. Walker, a  well-known differential geometer who had also reviewed Schouten's book of 1954 for {\it Mathematical Reviews}. There were three other papers written by Schouten jointly with Yano in the same year 1955 and those were reviewed by Nijenhuis.


\subsection*{The Schouten--Nijenhuis bracket in Poisson geometry}
Some twenty years elapsed before the Schouten--Nijenhuis bracket came to play a signifcant role in geometry. 
The history of the role that it had played in the infancy of Poisson geometry is well known\footnote{See Y. Kosmann-Schwarzbach, ed.,
{\it Sim\'eon-Denis Poisson : Les Math\'ematiques au
service de la science}, \'Editions de l'\'Ecole polytechnique, 2013, 147--168.}. Important parts were played by  
Andr\'e Lichnerowicz (1915--1998), as early as the Rome colloquium on symplectic geometry in 1973, and
by W\l odzimierz Tulczyjew in 1974 when he made use of ``the Lie bracket structure for multivectors'' to characterize those bivectors that are the inverse of a symplectic form. Then, in 1975, an article by  Lichnerowicz -- a study of a generalization of symplectic manifolds in terms of what he called ``le crochet de Nijenhuis'' -- appeared, followed by two papers written with Mosh\'e Flato (1937--1998) and Daniel Sternheimer to further study these generalized structures, and in 1976 a joint paper on ``Deformations of Poisson brackets, Dirac brackets and applications'' \cite{LFS1976} 
citing Nijenhuis for his 1955 paper and Gerstenhaber for his article on deformations of 1964. Eventually, in 1977, Lichnerowicz published a formal definition of Poisson manifolds in his article in the {\it Journal of Differential Geometry} where he proved that the vanishing of the Schouten--Nijenhuis bracket of a bivector is equivalent to the defining condition for a Poisson manifold, the Jacobi identity for the bracket of functions defined by the bivector.


\section{The Fr\"olicher--Nijenhuis bracket}
\subsection*{The Nijenhuis torsion and the bracket of $(1,1)$ tensors, 1951}
It was Nijenhuis who, in his article ``$X_{n-1}$-forming sets of eigenvectors'' \cite{N1951}, published in 1951, introduced the general expression for the torsion of a field 
of $(1,1)$-tensors -- then called mixed tensors of valence 2.
Schouten, in his 1951 paper ``Sur les tenseurs de $V^{n}$ aux directions principales $V^{n-1}$-normales'' (On the tensors of $V^n$ whose principal directions are $V^{n-1}$-normal)  \cite{Schouten1951}, had just studied the problem of the integrability of the eigendistributions spanned by $n-1$ fields of eigenvectors of a field of $(1,1)$-tensors on a manifold of dimension $n$, explaining and vastly extending a paper \cite{T1949} published two years earlier by the Italian geometer, Angelo Tonolo (1885--1962).
In \cite{N1951}, Nijenhuis in turn studied the integrability of the distributions spanned by the eigenvectors of a field of endomorphisms of the tangent spaces of a manifold.
He announced
``a hitherto unpublished differential comitant''
whose vanishing is ``equivalent to Tonolo--Schouten's criteria'' and permits expressing the conditions of integrability by the vanishing of quantities ``contain[ing] neither eigenvalues nor eigenvectors''.
He defined in local coordinates a quantity with components 
$H^{ . .  \kappa}_{\mu \lambda}$, expressed
in terms of the components $h^{. \kappa}_\lambda$ of a $(1,1)$-tensor, $h$, and of their partial derivatives, 
 $${  H^{ . \, .  \kappa}_{\mu \lambda} =  2 
 h^{~. \rho}_{[ \mu}   \partial^{}_{|\rho|}
 h^{. \kappa}_{\lambda]}  - 2 
  h^{. \kappa}_\rho \partial^{}_{[ \mu}
   h^{. \rho}_{\lambda ]}},$$
where the indices between square brackets are to be skew-symmetrized. He then proved the tensorial character of this quantity. He also introduced the bracket of a pair of $(1,1)$-tensor fields, $h$ and $k$,  
 $$[h,k]^{ . \, .  \kappa}_{\mu \lambda} = 
 h^{~. \rho}_{[ \mu}   \partial^{}_{|\rho|}
 k^{. \kappa}_{\lambda]}
+
 k^{~. \rho}_{[ \mu}   \partial^{}_{|\rho|}
 h^{. \kappa}_{\lambda]}  
   - 
  h^{. \kappa}_\rho \partial^{}_{[ \mu} 
   k^{. \rho}_{\lambda ]}
-
k^{. \kappa}_\rho \partial^{}_{[ \mu} 
   h^{. \rho}_{\lambda ]},
$$
the symmetric bilinear expression $[h,k]$ which yields the above quadratic form when $h = k$.
He proved that the bracket of a pair of $(1,1)$-tensor fields is a differential $2$-form with values in vector fields. 

Except for a normalising factor $2$, these were the expressions in local coordinates of the ``Nijenhuis torsion'' of a $(1,1)$-tensor field and the ``Nijenhuis bracket'' of a pair of $(1,1)$-tensor fields as they are now defined.
  
\subsection*{The Nijenhuis bracket of vector-valued differential forms, 1955}
Four years after he
had defined the torsion of a mixed tensor using local coordinates in \cite{N1951}, Nijenhuis formulated his previous results in coordinate-free fashion, defined the bracket of pairs of vector-valued differential forms of all degrees that generalized the bracket he had defined in 1951, and he proved that a graded Lie algebra structure was thus obtained on the space of vector-valued differential forms. This study  constitutes the latter part of his article \cite{N1955} which appeared in 1955,  -- the same paper in which he proved the Jacobi identity for the bracket of symmetric contravariant tensors and a graded Jacobi identity for the bracket of skew-symmetric contravariant tensors.

He introduced the concomitant $[h,k]$ of two $(1,1)$-tensors 
$ h^{\kappa}_{\cdot \lambda}$ and $k^{\kappa}_{\cdot \lambda}$
-- which he had defined in his 1951 paper -- first in local coordinates.  Then, referring to the formula for $[h,k]$ which he had  written there in terms of Lie derivatives, he wrote, in a coordinate-independent way, the formula for the bracket $[h,h]$, associated to a $(1,1)$-tensor $h$, evaluated on vectors $u$ and $v$,
$$
[h,h](u,v) = - h \cdot [h \cdot u,v] - h \cdot [u, h \cdot v] 
+ h \cdot h \cdot [u,v] + [h\cdot u, h \cdot v].$$
He pointed out that, for $h = k$, the concomitant $[h,h]$ is ``of special importance'' and, for that case, ``the name {\it torsion} seems to be quite appropriate''.
Giving an invariant form to a relation he had derived 
in \cite{N1951}, he wrote the identity
\begin{small}
$$
[h,k \circ l](u,v) + [k,h \circ l](u,v) = 
 k \cdot [h,l](u,v) + h \cdot [k,l](u,v) 
- [h,k](u, l \cdot v) - [h,k](l \cdot u, v),$$
\end{small} 
a ``formula [which] enables one to express $[h^q,h^q]$ algebraically in terms of $[h,h]$ and $h$''. Calling a field of $(1,1)$-tensors with vanishing torsion a ``Nijenhuis operator'', this formula, used recursively, proved that any power of a Nijenhuis operator is a Nijenhuis operator. 

In the last section of his article, ``A concomitant for differential forms with values in the tangent bundle'', he extended the definition of the bracket of $(1,1)$-tensors to vector-valued forms of arbitrary degrees, proving its tensorial character and stating the relevant ``Jacobi-type identity'',
$$(-1)^{rp}[L,[M,N]] + (-1)^{pq} [M,[N,L]] + (-1)^{qr}[N,[L,M]] =0,$$
for $L$ a vector-valued $p$-form, $M$ a vector-valued $q$-form, and $N$ a vector-valued $r$-form. 

He had remarked earlier, anticipating this general definition of a bracket for vector-valued forms of arbitrary degrees, that for any $(1,1)$-tensor, $h$, the equations   
$[h,[h,h]]=0$ and $[[h,h],[h,h]]=0$ are satisfied. 

\subsection*{Fr\"olicher and Nijenhuis, 1956} 
Alfred Fr\"olicher (1927--2010) was a 
Swiss mathematician who 
defended his thesis at the Eidgen\"ossische Technische Hochschule in Zurich in 1954 with a dissertation ``\selectlanguage{german}Zur Differential\-geometrie der komplexe Strukturen\selectlanguage{english}'' (On the differential geometry of complex structures) which was published in {\it Mathematische Annalen} the following year \cite{Frolicherthesis}. Later he was appointed to the Universit\'e de  Fribourg, then, in 1966, to the Universit\'e de Gen\`eve.

In 1956, Nijenhuis published another article \cite{FN1956}, written jointly with Fr\"olicher, on the properties of the bracket which he had introduced the previous year: the bracket of vector-valued differential forms. They characterized the two types of derivations 
of the exterior algebra of differential forms on a smooth manifold, proving that any derivation is the sum of a derivation that commutes with the exterior differential, defined by a vector-valued form (type $d_*$) and a derivation vanishing on functions (type $i_*$).
Whence the definition of a bracket of vector-valued forms: the bracket, $[L,M]_{\rm{FN}}$, of a vector-valued $\ell$-form, $L$, and a
vector-valued $m$-form, $M$, is the vector-valued $(\ell + m)$-form which satisfies the equation 
$${\mathcal L}_{[L,M]_{\rm{FN}}} = [{\mathcal L}_L, {\mathcal L}_M],
$$
where the bracket $[~,~]$ is the graded commutator of derivations of the algebra of differential forms, and ${\mathcal L}_U = [i_U,\d]$ is the graded commutator of $i_U$, the interior product by a vector-valued form $U$, and the de Rham differential,~$d$.
Thereafter, the bracket of vector-valued differential forms, which had been invented by Nijenhuis in 
, became known as the {\it Fr\"olicher--Nijenhuis bracket}.


\subsection*{Almost complex structures}
In his article of 1955 \cite{N1955}, Nijenhuis wrote, concerning the condition $[h,h] =~0$ that was necessary for a $(1,1)$-tensor field $h$ to define a complex structure, 
\begin{quotation}
\noindent ``This condition is also sufficient if the almost-complex structure is real analytic. For almost-complex structures of, say, class $C^\infty$, no satisfactory sufficient conditions for the structure to be a complex one, have yet been found,''
\end{quotation}
and in a footnote,
he observed:\begin{quotation}
\noindent It is rather remarkable that Eckmann and Fr\"olicher  
found $[h,h]$ in a slightly different form, so that it is a tensor only if $h \circ h = const$.\end{quotation}
Nijenhuis mentions that 
Yano had lectured on ``the relation between the two forms of the torsion'' in 1954, in his address at the
International Congress of Mathematicians 
which took place in Amsterdam.

Concerning the first comment, we now know that, in the $C^\infty$ case, the vanishing of the torsion is also sufficient and that the proof would be the celebrated Newlander--Nirenberg theorem which appeared in the {\it Annals of Mathematics} in 1957. When Nijenhuis was invited to give a 30-minute lecture at the next ICM in Edinburgh in 1958, he reviewed the recently proved theorem.

The comment in the footnote refers to the early history of the concept of torsion, which in fact appeared both 
in the work of the French mathematicians, Charles Ehresmann (1905--1979) and  Paulette Libermann (1919--2007), and in that of Beno Eckmann (1917--2008) together with 
Fr\"olicher, and in Fr\"olicher's thesis.

Eckmann, who had defended a brilliant thesis in 1941 at the Eidgen\"ossische Technische Hochschule in Zurich, first  taught at his {\it alma mater}, then at the Universit\'e de Lausanne. He visited Princeton in 1947 and returned to Switzerland to accept a chair of full professor at ETH in 1948.
In 1951, Eckmann together with his student   Fr\"olicher wrote a Note in the {\it Comptes rendus} of the Paris Academy of Sciences, ``Sur l'int\'egrabilit\'e  des structures presque complexes'' (On the integrability of almost complex structures) \cite{EF1951}. They introduced the torsion  in terms of local coordinates and
proved that these coordinates, $t^i_{kl}$, were those of a tensor 
and that the vanishing of this torsion is necessary for an almost complex structure to be a complex structure, and also sufficient in the analytic case. The technique of defining an object by means of an auxiliary connection and then proving that it was in fact independent of the choice of a connection was used frequently in differential geometry, as it still is today. Unlike the articles of Nijenhuis of 1951 and 1955, this paper only dealt with the case of a field of endomorphims of the tangent spaces to a manifold
whose square at each point is the opposite of the identity transformation.

The concept of the torsion of an almost complex structure had been introduced 
in 1950 by Ehresmann, professor at the University of Strasbourg, in his address to the International Congress of Mathematicians which was held at Harvard University, a study of the integrability of almost complex structures on manifolds. His doctoral student, Libermann, had already published two notes in the {\it Comptes rendus de 
l'Acad\'emie des Sciences} in 1949 with her adviser, another one followed in 1951. She published several short communications until her thesis -- which she had defended in 1953 -- was published in 1954 \cite{L1954}. 
That the torsion of an almost-complex structure, originally defined in terms of a connection, was in fact independent of the choice of a connection followed from both the approach of Libermann, in which the torsion was defined in 
terms of complex differential Pfaffian forms,
and that of Eckmann and Fr\"olicher  who worked ``staying in the framework of the real tensorial calculus'', as Libermann wrote in 1955 \cite{L1955}.
That the rival methods were announced and published nearly simultaneously was well explained in this article.

The expression of this torsion in terms of the  Lie brackets of vector fields was derived by Eckmann -- two years after his joint paper with Fr\"olicher -- in his communication, ``Sur les structures complexes et presque complexes'' (On complex and almost complex structures) \cite{BE1953} at the Colloque de G\'eom\'etrie diff\'erentielle held in Strasbourg in 1953.
We read in Chern's review  in {\it Mathematical Reviews} of Eckmann's paper that ``an equivalent condition is given, for an almost complex structure to be without torsion, which is in terms of vector fields and their commutator operators''. In fact, Eckmann's formula (3) is
$$\Omega(a,b) =  [a,b] + J [Ja, b] + J[a,Jb]  - [Ja,Jb], 
$$
free of the appearance of a connexion, and is now our modern version of the torsion, still in the particular case of $J^2 = -1$.
As for Libermann, at the same conference, she expressed the torsion in terms of the covariant differential 
with respect to a connexion, while Ehresmann chose to lecture on jet spaces, prolongation of structures, and Lie pseudo-groups, under the  title ``Introduction to the theory of infinitesimal structures and Lie pseudo-groups'', a title translated from the French since he spoke in French as did Chern, Kuiper, Lichnerowicz, Reeb, Souriau, Willmore and all the other distinguished speakers. 

When Fr\"olicher, who had attended the conference in Strasbourg,  published his doctoral thesis in 1955 \cite{Frolicherthesis}, he cited Eckmann, Ehresmann and Libermann. In the first chapters he used local coordinates to define the torsion of an almost complex structure, and, in a footnote, he cited the independent work of Nijenhuis who had defined a differential concomitant as a mixed tensor field that coincided with his own $t^j_{kl}$ in the case where $a^i_k a^k_j = - \delta^i_j$, i.e., when the mixed tensor has square $-1$. The expression of the torsion in terms of Lie brackets as above, still in the case where $J^2 = -1$, was given on page 77 of his thesis and used repeatedly.

It was the following year that Fr\"olicher collaborated with Nijenhuis in the joint paper \cite{FN1956} that was analysed above, dealing this time with the more general notion of torsion for all 
$(1,1)$-tensor fields that Nijenhuis had considered in 1951.

\section{Brackets and deformations}

\subsection*{The Gerstenhaber bracket, 1963}
One recent event has cast a bright light on the not so recent history of the Gerstenhaber bracket.
Murray Gerstenhaber was awarded the prestigious Steele prize of the American Mathematical Society on November 24, 2020, 
for the two remarkable articles of 1963 and 1964 which ``established the foundations of algebraic deformation theory''.

\m

Independently of the geometry in the context of which the Schouten--Nijenhuis bracket had appeared and eventually came to play such an important role, a new concept arose in algebra, more precisely in Gerstenhaber's theory of deformations of algebraic structures. That both theories dealt with a common object, graded Lie brackets, was immediately clear. 
In July 1962, Gerstenhaber submitted ``The cohomology structure of  an associative ring'' to the {\it Annals of Mathematics} and the article \cite{G1963} appeared the following year with the mention ``originally intended as the first section of a paper in preparation on the deformations of associative algebras''. In fact, his second paper \cite{G1964} was submitted soon after the first and was published early in 1964.

\m

In his  1963 paper \cite{G1963},
Gerstenhaber determined the structure of the cohomology ring $H^*(A, A)$ of a commutative ring $A$. The reviewer for {\it Mathematical Reviews}, Maurice Auslander (1926--1994), himself a distinguished algebraist, summarized the results as follows:
``Among other things, it is shown that $H^*(A, A)$ is a commutative ring in the sense of graded rings.
A bracket product $[\, , \,]$ is introduced in 
$H^*(A, A)$ under which 
$H^*(A, A)$  becomes a graded Lie ring with $[H^m(A,A),H^n(A,A)]\subset H^{m+n-1}(A,A)$ and such that the bracket operation of $H^1(A, A)$ into itself is the ordinary Poisson bracket [i.e., commutator] of derivations of $A$ into itself.''

Gerstenhaber first introduced {\it pre-Lie systems}, 
giving as an example the sequence of modules $V_m ={\mathrm{Hom}}(V^{\otimes (m+1)}, V)$, $m = -1, 0, 1, ...$, where $V$ is a module on a commutative ring, with the operations $\circ_i$ defined as follows. If $f$ and $g$ are of degree $m$ and $n$, respectively, in the sense that $f \in V_m$ and $g \in V_n$, then $f\circ_i g$ is an element in $V_{m+n}$ obtained by replacing the element in the $i$-th place, in the arguments of $f$, by $g$ evaluated
on the tensor product of the $n+1$ elements that follow the $i$-th argument. Summing over $i$ with the correct signs, that is all $+$ if the second factor is of even degree and an alternate sum otherwise, yields 
a {\it composition product} denoted by $f \circ g$ and he proved that this multiplication on $V = \oplus_m V_m$
 satisfies the ``graded right pre-Lie ring identity'', 
$$
(f \circ g) \circ h - f \circ (g \circ h) = (-1)^{np}
((f \circ h) \circ g -
f \circ (h \circ g)),
$$
for $f \in V_m$, $g \in V_n$, and $h \in V_p$.
Upon skew-symmetrizing the composition product in a graded pre-Lie ring, one obtains a graded Lie bracket,
henceforth denoted by $[~,~]$. 
Thus the module $C^*(A,A)$ of cochains on a module 
$A$ with values in $A$ is a Lie ring under 
the bracket $[~,~]$ which is, by definition,  the  skew-symmetrized composition product, $\circ$, itself
defined in terms of the operations $\circ_i$. 
Thus, when denoting ${\mathrm{Hom}}(A^{\otimes (m+1)}, A)$ by $A_m$, bracket $[~,~]$ maps a pair of cochains in $A_m$ and $A_n$, respectively, to a cochain in $A_{m+n}$.
 
When $A$ is an associative ring, the coboundary operator $\delta$ is defined on cochains on $A$ with values in $A$, whence the definition of the Hochschild cohomology $H^*(A,A)$ of the ring $A$. 
Since $A$ is an associative ring, the cup product 
$\smile$ is defined, on the space $C^*(A,A)$ of cochains on $A$ with values in $A$, by
$(f \smile g)(a,b) = f(a) g(b)$ when $f$ is an $m$-cochain, $g$ an $n$-cochain, $a \in A^{\otimes m}$ and  $b \in A^{\otimes n}$. 
On $H^*(A,A)$,  the cup product induces a product, denoted by the same symbol.
Thus, when equipped with the cup product, the cohomology of an associative ring is a graded ring called the ``cohomology ring of $A$''.
  
Still in the case of an associative ring, $A$, 
Gerstenhaber introduced a graded Lie bracket 
$[~,~]$, on $H^*(A,A)$, induced from the graded Lie bracket he had defined on $C^*(A,A)$.
He called $(H^*(A,A),[~,~])$ the ``infinitesimal ring of the ring $A$''.

Now, let $\pi : A \otimes A \to A$ denote the multiplication of the ring $A$ considered as a $2$-cochain on $A$ with values in $A$. Because the multiplication is associative, $\pi$ is a cocycle, ``which may be called the canonical $2$-cocycle of $A$''. Gerstenhaber proved two facts concerning $\pi$.
On the one hand, for cochains $f$ and $g$ on the ring $A$,
$f \smile g = (\pi \circ f) \circ g$, a relation expressing the cup product of cochains in terms of the multiplication $\pi$, considered as a $2$-cochain,  and the composition $\circ$.
On the other hand, the coboundary operator $\delta$  is, up to sign, the right adjoint action of $\pi$ with respect to the bracket,
$\delta f =  - [f,\pi]$, for any cochain $f$ on $A$, and in particular, $[\pi, \pi] =0$.
Thus, the cohomology operator, $\delta$, is ``a right inner derivation of degree $1$ of the graded Lie ring $(C^*(A,A), [ ~,~])$''. 

Then, Gerstenhaber proved that the defect in the derivation property of
the coboundary operator $\delta$ with respect to the $\circ$ multiplication is given by the cup bracket defined as the skew-symmetrized cup product.
As a consequence, the ring structure, $\smile$, on $H^*(A,A)$ induced by the cup product  is graded commutative.
From the space of cochains of an associative ring $A$, the Hochschild cohomology, $H^*(A,A)$, inherits the structure of a graded Lie ring under the 
bracket, $[~,~]$. Finally, the interplay of the cup product and the bracket is demonstrated: for cocyles, $f, g, h$, of degrees $m,n,p$, respectively,
$[f \smile g,h] - [f ,h] \smile g - (-1)^{m(p-1)} f \smile [g,h]$ 
is a coboundary, whence the vanishing of this expression in  $H^*(A,A)$ which expresses the fact that, for $\xi \in H^p(A,A)$, the map $\eta \to [\eta, \xi]$ is a left derivation of degree $p-1$ of the cup product.
It is in Corollary 2 of Theorem 5 that these facts are proved and summarized in the statement that for any $\xi \in H^*(A,A)$, $D_\xi : \eta \mapsto [\eta,\xi]$, is a left derivation of $(H^*(A,A), \smile)$ and $\xi \mapsto D_\xi$ is an anti-homomorphism of rings. Here one can recognize a ``Gerstenhaber algebra structure'' on $(H^*(A,A), [~,~], \smile)$.
  
\m

In fact, Gerstenhaber dealt more generally with the cohomology module $H^*(A,P)$, where $P$ is a two-sided $A$-module, and proved that,
$H^*(A,P)$ is a two-sided module over the Lie ring 
$(H^*(A,A),[~,~])$.

\m

He gave a useful summary of his results as an introduction to his article, collecting its main formulas: if $A$ is an associative ring, the module of cochains on $A$ with values in $A$ is then equipped with both the cup product and the composition bracket which is shown to be skew-commutative, ``the grading being reduced by one from the usual''. The cohomology $H^*(A,A)$ inherits these structures from the module of cochains, the cup product $\smile$ which becomes graded commutative and the bracket, $[~,~]$, induced from the composition bracket of cochains, which is skew-symmetric  with respect to the shifted grading, yielding, for cohomology classes, $\alpha$, $\beta$ and $\gamma$,
$$
[\alpha, \beta] = - (-1)^{(m-1)(n-1)} [\beta,\alpha]
$$
and is shown to satisfy the following graded Jacobi identity,
$$
(-1)^{(m-1)(p-1)} [[\alpha,\beta],\gamma]
+ (-1)^{(n-1)(m-1)} [[\beta,\gamma], \alpha]
+ (-1)^{(p-1)(n-1)} [[\gamma,\alpha],\beta]=0,
$$
where $m,n, p$ are the degrees of  $\alpha, \beta, \gamma$,  respectively.
The restriction to elements of degree $1$ of this ``bracket product'' is the commutator of derivations.
 In addition, the bracket and the cup product were shown to be related by the identity 
$$
[\alpha \smile \beta, \gamma] = 
[\alpha , \gamma]\smile \beta + (-1)^{m(p-1)} \alpha \smile [\beta, \gamma],
$$
for $\alpha$  of degree $m$ and $\gamma$ of degree $p$,
which showed that the map $\alpha \mapsto [\alpha, \gamma]$ is ``a derivation of degree $p-1$ of $H^*(A,A)$ considered as a ring under the cup product''.
To summarize further, Gerstenhaber showed that the cohomology $H^*(A,A)$ of an algebra $A$ with values 
in $A$ is a graded commutative algebra with respect to the cup product multiplication $\smile$ and a graded Lie algebra with shifted grading with respect to the bracket, and that, in addition, the two structures are related by the identity that characterizes a ``Gerstenhaber algebra''.

\m

Years later, in 1988, in ``Algebraic cohomology and deformation theory'' \cite{GS1988}, a survey of past and recent results  which he published with his former doctoral student and collaborator Samuel D. Schack (1953-2010), 
Gerstenhaber recalled that ``algebraic deformation theory `broke free' in the spring of 1961'' when Andr\'e Weil conjectured the existence of a structure yet to be discovered on the Hochschild cohomology of an algebra, and that 
``the concept of `graded Lie algebra' had not previously been formalized but had obvious significance and was thus published before the deformation theory [i.e., in \cite{G1963}, the first of  Gerstenhaber's two papers]. Graded Lie algebras were then recognized everywhere.''
It was not until 1992 that the expression ``Gerstenhaber algebra'' appeared on page 8 of the article ``Bimodules, ...'' that Schack published in the {\it Journal of Pure and Applied Algebra}\footnote{This reference comes from the unpublished thesis of Kristin Haring, ``On the Events Leading to the Formulation of
the Gerstenhaber Algebra: 1945--1966'' (University of North Carolina, Chapel Hill, 1995),
which contains a wealth of additional information.}. Then, Bing~H.~Lian and Gregg Zuckerman, in their BRST paper in {\it Commun. Math. Phys.}, received in December 1992 and published in 1993, wrote of ``what  we [the authors] call  the  Gerstenhaber  bracket'' and ``what  mathematicians  call  a  Gerstenhaber  algebra'' in their abstract. That is when
the term ``Gerstenhaber algebra'' became standard.

\m

The left pre-Lie algebras were later variously called left symmetric algebras or Vinberg algebras. However, the name ``pre-Lie algebras'' remains current. 
In Gerstenhaber's 1963 paper, we find the definition of a bracket on the module of cohomology classes of cochains on a module over a ring. When the module is the module of smooth vector fields over the ring of smooth functions on a smooth manifold, one recovers the Schouten--Nijenhuis bracket of multivector fields. 

\m 

Gerstenhaber's article of 1964 \cite{G1964} developed a ``deformation theory for rings and algebras'',  considering mainly associative rings and algebras, ``with only brief allusions to the Lie case, but the definitions hold for wider classes of algebras.''
He recalled that
``certain aspects of the present deformation theory parallel closely those
 of the Fr\"olicher--Kodaira--Nijenhuis--Spencer theory'',
 and he spelt out the general principle of 
 ``identification of the infinitesimal deformations of a given object
 with the elements of a suitable cohomology group.''
 He concluded his introduction stressing the significance of his first paper since ``the deformation theory for algebras has shown that the direct sum of the
 groups $H^n(A, A)$ possesses a much richer structure than had previously
 been exhibited'', in fact a ``Lie product'' on $H^*(A,A)$, i.e., a Gerstenhaber algebra structure on the Hochschild  cohomology of a ring. He wrote that ``an analogous operation is definable for Lie rings,'' citing the still unpublished work of Nijenhuis and Richardson.
He summarized the three chapters that follow. In the first, the interpretation of 
$H^2(A, A)$ as the group of infinitesimal
deformations of $A$ is given ``in precisely the same way that the first cohomology group (derivations of A into itself modulo inner derivations) is interpreted as the group of infinitesimal automorphisms.'' In the second chapter,
he discussed ``the set of structure constants for associative algebras of dimension $n$ as a parameter space for the deformation theory of these algebras'', while the last chapter developed a
deformation theory for graded and filtered rings with applications to the rigidity of algebras.


\m

In the conclusion of \cite{GS1988}, Gerstenhaber  recalled that, following his 1963 paper on deformations of associative algebras, there remained domains of deformation theory to be developed. He wrote:
\begin{quotation}
``That deformation theory was meaningful for other categories [i.e., other than the associative algebras], e.g. Lie algebras, was explicit in \cite{G1963} -- and, in any case obvious. The challenge to develop it in the Lie case was rapidly picked up by Nijenhuis and Richardson.''
\end{quotation}
I shall now describe how these authors dealt with this challenge.


\subsection*{Nijenhuis and Richardson, 1964}
Nine years after the 1955 paper \cite{N1955} which introduced the graded Lie algebra of vector-valued forms, Nijenhuis with Roger W. Richardson (1930-1993) developed the theory of a new graded bracket.
Richardson  was an algebraist, who had defended his Ph.D. on group theory at the University of Michigan, Ann Arbor, in 1958. In 1966, he was teaching at the University of Washington in Seattle, when 
Nijenhuis was already at the University of Pennsylvania in Philadelphia but we know that he visited Seattle often.
In 1964, they published a research announcement, ``Cohomology and deformations of algebraic structures'', in the {\it Bulletin of the American Mathematical  Society} \cite{NR1964}, and two years later a full-length article appeared in the same journal, ``Cohomology and deformations in graded Lie algebras'' \cite{NR1966} -- based on Nijenhuis's address, ``Derivations and structures'',  
at the 1962 meeting of the American Mathematical Society in Vancouver --, followed by an expository article in 1967 in the {\it Journal of Mathematics and Mechanics}  \cite{NR1967}, which was ``a delight to read'' as the reviewer for {\it Mathematical Reviews} described it.
Unlike the Schouten--Nijenhuis and Fr\"olicher--Nijenhuis  brackets which were chapters in the differential geometry of manifolds, 
the  bracket of Nijenhuis and Richardson was an algebraic construct motivated by the theory of deformations of algebras.

In their short paper of 1964 \cite{NR1964}, which they had submitted at the end of 1963, there is no abstract and the paper opens with: 
\begin{quotation}
\noindent ``Gerstenhaber has recently initiated a theory of deformations of associative algebras. The  methods  and  results  of  Gerstenhaber's  work  are  strikingly  similar  to  those  in  the  theory  of  deformations  of  complex  analytic  structures  on  compact  manifolds.''
\end{quotation}
They announced that, in this note, they would reformulate Gerstenhaber's ideas within the framework  of  graded  Lie  algebras. Citing the article of Kunihiko Kodaira, Louis Nirenberg and Donald C. Spencer of 1958 and the 1962 article by Masatake Kuranishi, both in the {\it Annals of Mathematics},  on the existence of deformations of complex analytic structures, they distinguished the ``basic role played by a certain equation among the $1$-dimensional elements of a graded Lie algebra [which] expresses the integrability conditions for almost complex structures'', and they stated:
``Our basic observation  is that  a wide class  of algebraic  structures  on  a  vector  space  can  be  defined  by  essentially  the  same  equation  in  an  appropriate  graded  Lie  algebra''. 
They proceeded to describe a general theory of deformations of algebraic structures and applications.
Their aim was then to study the rigidity of structures, obtaining for Lie algebra structures an analogue of Kuranishi's result of 1962 for complex structures, as well as an analogue of the recently published theorem of Gerstenhaber for associative algebras  \cite{G1964}. 


Nijenhuis and Richardson first defined graded Lie algebras, and they stated that, for several types of algebraic structures on a given vector space,~$V$, there is a graded Lie algebra $(E, [~,~])$ associated to $V$ such that an element, $x$, of degree $1$ in $E$ satisfying the equation 
$$
[x,x] =0
$$
is an
algebra structure on $V$ of the given type, and then the linear map $y \mapsto [x,y]$ of $E$ into $E$, denoted by $\delta_x$, is a derivation of degree $1$ and square~$0$ of $E$, and the direct sum, $H(x)$, for all non-negative degrees, of the cohomology spaces defined by $\delta_x$  ``inherits from $E$ the structure of a graded Lie algebra''. 

As an example of an algebraic structure thus defined in an appropriate graded Lie algebra, they treated the case of Lie algebras.
Given a vector space, $V$, on a field of characteristic
$\neq 2$, they defined the hook product on $A(V)= \oplus_{n\geq - 1} \, A^n(V)$, where $A^n(V)$ is the space of skew-symmetric
multilinear maps from $V^{n+1}$ to $V$. The {\it hook product}  of $f \in A^p(V)$ and $g \in A^q(V)$ is the element $f \barwedge h$ in  $A^{p+q}(V)$ whose evaluation on $(u_0, u_1, \dots, u_{p+q})$ is the sum, over all permutations $\eta$ of 
$\{0, 1, ..., p+q\}$ that map both ${0, 1,  ..., p}$  and ${p+1, p+2, ..., p+q}$ to strictly increasing sequences,
of the terms ${\mathrm{sgn}} (\eta) f(h(u_{\eta(0)}, ...,u_{\eta(p)}), u_{\eta(p+1)}, ... , u_{\eta(p+q)})$ where ${\mathrm{sgn}} (\eta)$
is the signature of the permutation $\eta$.
Then the {\it hook bracket} is defined as the skew-symmetrized hook product and it is a graded Lie bracket on $A(V)$, denoted by $[~,~]$. 
Thus, the hook bracket, nowadays called the Nijenhuis--Richardson bracket, is a graded Lie bracket 
of degree $-1$ on the space of vector-valued forms on any vector space on a field of characteristic $\neq 2$.

With this bracket, they characterized the Lie algebra structures on a vector space, $V$, as those elements $f \in A^1(V)$, i.e., linear maps from $\wedge^2 V$ to $V$,
that satisfy the condition $[f,f] = 0$, and they considered, for a Lie algebra structure $f$ on $V$, the  cohomology operator $\delta_f$, acting on the space of vector-valued forms on $V$ as $\delta_f = [f,  \cdot]$. Then, they stated that the cohomology $H(f)$ coincides, up to a shift of $1$ in degree, with the Chevalley--Eilenberg cohomology of the Lie algebra $L= (V,f)$, 
$$
H^n(f) = H^{n+1}(L,L).
$$
They stated that ``similar remarks hold for the case of associative multiplications on $V$'', for which they cited Hochschild's paper of 1945 \cite{Hochschild1945}, as well as for associative and commutative multiplications for which they referred to David K. Harrison's 1962 article   \cite{Harrison1962}. 
They wrote: ``The existence of a graded Lie algebra structure on $H(L,L)$ ($= \sum H^n(L,L)$)
is due to Gerstenhaber \cite{G1963}.''


The fully developed account of their 
theory that Nijenhuis and Richardson gave two years later in \cite{NR1966} contained first a comprehensive survey of graded Lie algebras and then ``a detailed discussion of the {\it deformation
equation} in graded Lie algebras''.  
When $(E, [~,~])$ is a graded Lie algebra where each summand $E^n$ is finite-dimensional, given a derivation, $D$, of degree $1$ on $E$ such that $D^2=0$, the {\it deformation equation} is the equation
$$
Da + \frac{1}{2} [a,a]=0
$$
for elements $a \in E^1$. The authors stated that 
``in specific applications'', the set of solutions of this equation will be the set of structures of a given type on a given vector space and they illustrated their general statement. They pointed out that the special case where $D=0$ is important in applications. They discussed the set $M$ of solutions $a$ of the deformation equation which coincides with the set $N$ 
of solutions $x$ of 
the equation $[x,x]=0$ for $x=D+a$. 
They concluded that, if $m$ is in $N$, then $m+u$ is again in $N$ if and only if $u$ satisfies the equation 
$D_m u + \frac{1}{2} [u,u] =0$, where $D_m = D + [m, .]$, so that ``the `small' solutions of this deformation equation relative to $m$ give a neighborhood of $m$ on $M$''. Then they recalled the results on structures of algebras that were contained in their announcement of 1964 \cite{NR1964}, together with the relevant references \cite{Hochschild1945} \cite{Harrison1962} \cite{G1963}  \cite{G1964}, as well as the 1948 paper of Chevalley and Eilenberg, and they indicated an application of the deformation equation to the study of homomorphims of algebras. In the following sections, they treated the case of analytic, then algebraic graded Lie algebras, including rigidity theorems.

The ``deformation equation'' is what came to be known as the ``Maurer-Cartan equation'' which had been written, in terms of infinitesimal transformations, by Ludwig Maurer (1859--1927) in 1892, and ten years later, in terms of the exterior derivatives of forms, by \'Elie Cartan.


\section*{Conclusion}
The history of brackets did not end with the recognition that the Schouten--Nijenhuis bracket of 
differential geometry is a geometric counterpart of Gerstenhaber's bracket in deformation theory.
In the line of research originating in Schouten's ``concomitants'', which evolved into ``brackets'', new developments continued to appear, first in an article by Nijenhuis in 1967, 
and then in a paper that he wrote with Richardson in 1968,
followed by yet another paper that he published in 1969, ``On a class of common properties of different types of algebra'' 
in which his previous results were vastly generalized.
I~recall meeting Nijenhuis much later, in 1991, in Seattle -- at a large American Mathematical Society conference organized by, among others, Jerry Marsden --, and 
corresponding with him.
In a letter dated 1 May, 1998, reflecting on the Schouten bracket, he wrote: 
\begin{quote}
``After all, it is just an extension of the Lie bracket that, with some hindsight, anyone could have found. By contrast, the \break F-N
[Fr\"olicher--Nijenhuis] bracket is a bit more sophisticated, and its applications are fewer. -- Of course, I am keeping an eye on attempts to collect all these invariants under one heading,''
\end{quote}
and he added that ``in some of his last papers, Schouten also tried such a thing, not with much success''. As a matter of fact, a sort of ``bracket factory'' was developing at the time with the work of many younger mathematicians. I shall cite only Alexander Vinogradov and Theodore (Ted) Voronov, and I shall add that, by exploiting the concept of a Leibniz algebra defined by Jean-Louis Loday in 1993 and introducing the concept of ``derived bracket'' in 1995 -- after a suggestion by Jean-Louis Koszul contained in notes that he had written in 1990 and sent to me in 1994 --,
I was able to show how the various brackets that Schouten, then Nijenhuis, then Fr\"olicher and Nijenhuis, and then Nijenhuis and Richardson had defined were related\footnote{See
``Exact Gerstenhaber algebras and Lie bialgebroids'', {\it  Acta.
Appl. Math.} 41 (1995), and ``Derived brackets'', 
{\it Lett. Math. Phys.} 69 (2004).}. Voronov's theory of higher derived brackets was the next step in what had begun as a search for coordinate-independent differential expressions involving tensors.

\subsection*{Mackenzie's legacy}
Kirill Mackenzie's contributions to the development of the theory of Lie algebroids and Lie groupoids began with his four-page article in 1978 on the ``Rigid cohomology of topological groupoids'' \cite{M1978}
and continued with his 1979 thesis entitled {\it Cohomology of Locally Trivial Lie Groupoids and Lie Algebroids} which has remained unpublished.
When perusing his 1987 book, {\it Lie Groupoids and Lie Algebroids in Differential Geometry} \cite{Mackenzie1987}, we find no mention of graded brackets, no mention of Schouten, nor of Gerstenhaber, nor of any of the other authors whose names have become attached to these objects in differential geometry and homological algebra. The development of Poisson geometry since the mid-1970's and the theory of Poisson groupoids after they were first defined in 1988 by Alan Weinstein, created the need for algebro-geometric constructions of various kinds of brackets, especially those which had first been envisioned by Schouten and then fully developped by Nijenhuis. 
As early as 1995, Mackenzie had written about the need to introduce ``certain graded algebraic structures related
to Lie pseudoalgebras and Lie algebroids''\cite{1995}. 
So, when he was drafting his {\it General Theory of Lie Groupoids and Lie Algebroids}
\cite{Mackenzie2005},  it was clear to him that ``double groupoids arise naturally in Poisson geometry'', whence the need for some of the basic features of this geometry -- viewed from the point of view of the graded bracket structure of the algebra of  fields of multivectors -- in the construction of what Voronov called ``Mackenzie theory''\footnote{See Voronov's ``$Q$-manifolds and Mackenzie theory'', {\it  Comm. Math. Phys.} 315 (2012), and his luminous review of Mackenzie's {\it General Theory...} in the {\it Bulletin} of the London Mathematical Society of 2010.}.  No account of Poisson groupoids or of the infinitesimals of symplectic groupoids could be set out without recourse to those already classical ``concomitants'' in their modern guise, in particular the Poisson structures and their graded versions.

What Kirill Mackenzie achieved was to translate theorems, some old and many new, into the elegant language of morphisms between suitably defined objects.
His methods presented us with novel and powerful ways of dealing with those brackets, replacing many of their properties by conditions on the commutativity of beautiful, sometimes intricate diagrams.  By this method he was able to advance the theory of higher structures in vector bundle theory. In a sense, he taught us how to deal with the properties of brackets in terms of diagrams.

\noindent {\large{\bf Appendix.}} The list of references in the previous pages is limited to those items that have been explicitly cited in the text. 
The extent of the relevant literature since the 1970's is enormous. Below, I cite a small sample of important later papers, in abbreviated references, in chronological order.

C. Buttin, Théorie des opérateurs différentiels gradués sur les formes différentielles, {\it Bull. Soc. Math. France} 102 (1974) \begin{small}[article written by Pierre Molino ``based on ideas of Claudette Buttin'' (1935--1972)]\end{small}.

J.-L. Koszul, Crochet de Schouten-Nijenhuis et cohomologie, {\it Ast\'erisque} 1985.

I. S. Krasil'shchik,
Schouten bracket and canonical algebras,
{\it Lecture Notes in Math.} 1334 (1988).

A. M. Vinogradov, The union of the Schouten and Nijenhuis brackets, cohomology, and superdifferential operators (Russian), {\it Mat. Zametki} 47 (1990).

P. Lecomte, P. Michor, H. Schicketanz, The multigraded Nijenhuis-Richardson algebra, its universal property and applications,
{\it J. Pure Appl. Algebra} (1992).

J. Stasheff, The intrinsic bracket on the deformation complex of an associative algebra, {\it J. Pure Appl. Algebra} (1993).

J. Huebschmann, Lie-Rinehart algebras, Gerstenhaber algebras and Batalin-Vilkovisky algebras, {\it Ann. Inst. Fourier} (1998).

A. Cattaneo, D. Fiorenza, R. Longoni, Graded
Poisson algebras, {\it Encycl. Math. Physics} (2006).

A. Giaquinto, Topics in algebraic deformation theory,  {\it Prog. in Math.} 
287 (2011).

\bigskip

\noindent {\large{\bf Acknowledgments.}}
The judicious advice of Murray Gerstenhaber led to a definite improvement on the first version of this text. I am deeply thankful to him. 
To Jim Stasheff,
I express many thanks for his detailed, favorable comments and his always useful  questions. 
I thank Ted Voronov who responded to my first draft with encouragements and suggestions, while previous correspondence on related matters with Joseph Krasil'shchik,
Tudor Ratiu and Alan Weinstein is also gratefully acknowledged.

\flushright yks@math.cnrs.fr

\begin{thebibliography}{40}



\bibitem{CNS1975} L. Corwin, Y. Ne'eman, and S. Sternberg,
Graded Lie algebras in mathematics and physics (Bose-Fermi symmetry),
{\it Rev. Modern Phys.} 47 (1975), 573--603.  

\bibitem{BE1953} B. Eckmann, Sur les structures complexes et presque complexes, in  {\it G\'eom\'etrie diff\'erentielle} (Colloques Internationaux du Centre National de la Recherche Scientifique, Strasbourg, 1953), Centre National de la Recherche Scientifique, Paris, 1953, 151--159. 


\bibitem{EF1951}  B. Eckmann and A. Fr\"olicher, Sur l'int\'egrabilit\'e des structures presque complexes, {\it  C. R. Acad. Sci. Paris} 232 (1951), 2284--2286.
 
\bibitem{LFS1976} M. Flato, A. Lichnerowicz, and D. Sternheimer,
Deformations of Poisson brackets, Dirac brackets and applications,
{\it J. Math. Phys.} 17 (9) (1976),  1754--1762.

\bibitem{Frolicherthesis} A. Fr\"olicher, Zur Differentialgeometrie der komplexen Strukturen, {\it Math. Ann.} 129 (1955), 50--95.
 
\bibitem{FN1956} A. Fr\"olicher and A. Nijenhuis, Theory of vector-valued differential forms. I. Derivations of the graded ring of differential forms, {\it Nederl. Akad. Wetensch. Proc.} Ser. A. 59 = {\it Indag. Math.} 18 (1956), 338--359.

\bibitem{G1963}
M. Gerstenhaber, The cohomology structure of  an associative ring, {\it Annals of Mathematics}, Second Series, 78 (1963), 267--288.


\bibitem{G1964}
M. Gerstenhaber, On the deformation of rings and algebras, {\it Annals of Mathematics}, Second Series, 79  (1964), 59--103.

\bibitem{GS1988}
M. Gerstenhaber and S. D. Schack, Algebraic cohomology and deformation theory, in {\it  Deformation Theory of Algebras and Structures and Applications}, NATO Adv. Sci. Inst. Ser. C Math. Phys. Sci., 247, M.~Hazewinkel and M. Gerstenhaber, eds., Kluwer Acad. Publ., Dordrecht, 1988, 11--264. 


\bibitem{Harrison1962} D. K. Harrison, Commutative algebras and cohomology, {\it Trans. Amer. Math. Soc.} 104 (1962), 191--204. 


\bibitem{Hochschild1945} G. Hochschild, 
On the cohomology groups of an associative algebra.
{\it Ann. of Math.} (2) 46 (1945), 58--67.

\bibitem{L1954} P. Libermann, 
Sur le probl\`eme d'\'equivalence de certaines structures infinit\'esimales, 
{\it Ann. Mat. Pura Appl.} (4) 36 (1954), 27--120.

\bibitem{L1955} P. Libermann, 
Sur les structures presque complexes et autres
structures infinit\'esimales r\'eguli\`eres, {\it 
Bull. Soc. Math. France} 83 (1955), 195--224.

\bibitem{M1978} K. Mackenzie,
Rigid cohomology of topological groupoids,
{\it J. Austral. Math. Soc.} Ser. A.  26 (1978), no. 3, 277--301.

\bibitem{Mackenzie1987} K. C. H. Mackenzie,
{\it Lie Groupoids and Lie Algebroids in Differential Geometry}, London Mathematical Society Lecture Note Series 124, Cambridge University Press, Cambridge, 1987.

\bibitem{1995} K. C. H. Mackenzie, Lie algebroids and Lie pseudoalgebras,
{\it Bull. London Math. Soc.} 27 (1995), 97--147.
 
\bibitem{Mackenzie2005} K. C. H. Mackenzie,
{\it General Theory of Lie Groupoids and Lie Algebroids}, London Mathematical Society Lecture Note Series 213, Cambridge University Press, Cambridge, 2005. 


\bibitem{N1951} A. Nijenhuis, $X_{n-1}$-forming sets of eigenvectors, {\it Nederl. Akad. Wetensch. Proc.} Ser. A. 54 = {\it Indag. Math.} 13 (1951), 200--212.

\bibitem{N1952} A. Nijenhuis, {\it Theory of the geometric object}, Thesis, University of Amsterdam, Amsterdam, 1952.

\bibitem{N1955} A. Nijenhuis,
Jacobi-type identities for bilinear differential concomitants of certain tensor fields. I, II,
{\it Nederl. Akad. Wetensch. Proc.} Ser. A. 58 = {\it Indag. Math.} 17 (1955), 390--397, 398--403.  

\bibitem{N1972} A. Nijenhuis, J. A. Schouten: a master at tensors (28 August 1883--20 January 1971), {\it Nieuwe Arch. Wisk.} (3) 20 (1972), 1--19. 

\bibitem{NR1964} A. Nijenhuis and R. W. Richardson, 
Cohomology and deformations of
algebraic structures, {\it Bull. Amer. Math. Soc.} 70 (1964), 406--411. 

\bibitem{NR1966}  A. Nijenhuis and R. W. Richardson,  Cohomology and deformations in graded Lie algebras, {\it Bull. Amer. Math. Soc.} 72 (1966), 1--29. 

\bibitem{NR1967} A. Nijenhuis and R. W. Richardson,  Deformations of Lie algebra structures, {\it J. Math. Mech.} 17 (1967), 89--105.

\bibitem{N1918} E. Noether,
Invarianten beliebiger Differentialausdr\"ucke, {\it G\"ott. Nachr.} (1918), 37--44.
 
\bibitem{Schouten1914} J. A. Schouten, 
{\it Grundlagen der Vektor- und Affinoranalysis},
Leipzig, B.~G.~Teubner, 1914.

\bibitem{Schouten1940}  J. A. Schouten, Ueber Differentialkomitanten zweier kontravarianter Gr\"ossen,  {\it Nederl. Akad. Wetensch. Proc.} 43 (1940) = {\it Indagationes} 2 (1940), 182--185.

\bibitem{Schouten1951}  J. A. Schouten, 
Sur les tenseurs de $V^{n}$ aux directions principales $V^{n-1}$-normales, 
in {\it Colloque G\'eom. diff., Louvain 1951}, Centre Belge Rech. math. (1951), 67--70.

\bibitem{Schouten2} J. A. Schouten, {\it Ricci-Calculus: An Introduction to Tensor Analysis and its Geometrical Applications}, Springer, 1954. 

\bibitem{Schouten1954} J. A. Schouten, On the differential operators of first order in tensor calculus, in {\it Convegno Internaz. Geometria Differenz., Italia, 20--26 Settembre 1953}, Roma, Edizioni Cremonese (1954), 1--7 (Report 1953-012, Mathematische Centrum, Amsterdam).

\bibitem{SchoutenStruik} J. A. Schouten and D. J. Struik, {\it Einf\"uhrung in die neueren Methoden der Differentialgeometrie} (vol.~1, Algebra und \"Ubertragungslehre, von J.~A.~Schouten; vol.~2, Geometrie, von D. J. Struik), Groningen-Batavia, P.~Noordhoff N. V., 1935--1938. 

\bibitem{SchoutenYano1955b} J. A. Schouten and K. Yano, On the geometrical meaning of the vanishing of the Nijenhuis tensor in an $X_{2n}$ with an almost complex structure, {\it Nederl. Akad. Wetensch. Proc.} Ser. A. 58 = {\it Indag. Math.} 17 (1955), 133--138.
 
\bibitem{WS1931} W. {\'S}lebodzi{\'n}ski,  Sur les \'equations de Hamilton, {\it Bull. Acad. Royale de Belgique}
17 (1931), 864--870.

\bibitem{S1978}\noindent  D. Struik, Schouten and the tensor calculus, {\it Nieuwe Arch. Wisk.} (3) 26 (1978), 96--107.

 
\bibitem{T1949} A. Tonolo, Sopra una classe di deformazioni finite,
{\it Ann. Mat. Pura Appl.}, IV. Ser., 29 (1949), 99--114.

\end{thebibliography}
\end{document}